# Constructive representation of primitive Pythagorean triples

*Aleshkevich Natalia V.*

*Peter the Great St. Petersburg Polytechnic University*

**Abstract**

The paper presents a systematic construction of primitive Pythagorean triples. The order of enumeration on the set of primitive Pythagorean triples is defined. The order is based on the representation of a primitive Pythagorean triple by dividing the side of the generating square into two groups of factors using gnomons.

In paper is shown the inverse mapping of the elements of a primitive Pythagorean triple to the parameters of the partition of the side of the generating square, the representation of primitive Pythagorean triples as the sum of two connected gnomons, a one-to-one correspondence of the connected gnomons of the primitive Pythagorean triple with the corresponding arithmetic progressions. Also the property of such arithmetic progressions as connectedness with each other, namely, partial overlap, is found. The motion (configuration change) of the connected gnomons is shown during the transition of a primitive Pythagorean triple to a general Pythagorean triple, all the terms of which have one common coefficient.

**Setting the order on a set of primitive Pythagorean triples**

The construction of primitive Pythagorean triples using gnomons is shown in [1]. The entire set of primitive Pythagorean triples can be constructed according to the sequential growth of the known parameters. These parameters are the side of the generating square **S** and two groups of factors **t** and **l**. Since the parameters **t** and **l** are mutually related, we will choose one of them for ordering, namely **t**. Thus, the order is set by two parameters. One parameter external is the side of the generating square. Side **S** is even number. **S** starts with 2 and goes in increments of 2. The internal parameter **t** is the partition element of the side of the generating square **S=2tl**. The element **t** starts with the minimum value corresponding to the parameter **S**, and then increases to the maximum value within **S**. The associated element **l**, starting from the maximum, decreases at the same time. Both elements are formed from the cofactors of the number **S**. The number **t** can be of any parity. The number **l** is odd. Wherein $GCD(t, l) = 1$.

Formulas for obtaining elements of a primitive Pythagorean triple:

$$y = S + 2t^2;$$

$$x = S + l^2;$$

$$z = S + 2t^2 + l^2.$$

According to the parameters **S, t(S)** ordered tables of primitive Pythagorean triples can be constructed. In this case, the **N** ordinal number of the first level is equal to

$$N = S/2.$$

The sequence number **n** of the second level changes from 1 to **k** within **S**, where

$$k = \sum_{i=0}^{j} C_j^i = 2^j,$$

**j** is the number of prime odd cofactors without taking into account their powers included in the product for **S**.

The table is presented in Appendix 1. A fragment of the beginning of the set is given for the values $S = 2 \div 100$.

Accordingly $N = 1, 2, \ldots, 50.$

Setting the order allows to build different algorithms when using primitive Pythagorean triples.

**Mapping the elements of a primitive Pythagorean triple through the parameters of the partition of the side of the generating square**

We are given a primitive Pythagorean triple

$$x^2 + y^2 = z^2, \ y - even$$

Let's find the side of the generating square **S** and the elements of its partition **t, l**.

The partition parameter **l** is obtained as the square root of the difference between the hypotenuse and the even leg:

$$z - y = l^2, \quad l = \sqrt{z - y};$$

The difference between the odd leg and the square **l²** gives the value of the side of the generating square **S**:

$$x - l^2 = 2lt = S,$$

Next, we immediately get the second element of the partition **t**:

$$t = \frac{S}{2l}.$$

An example of mapping the elements of a primitive Pythagorean triple through the parameters of the partition of the side of the generating square is considered in Appendix 2.

**Primitive Pythagorean triples and their geometric representation**

Each square can be represented as a smaller square and a gnomon placed on this smaller square (Figure 1).

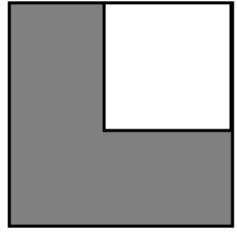

Figure 1

A primitive Pythagorean triple can be represented as a square with an even side and a gnomon whose area is equal to the area of a square with an odd side (Figure 2.1). And vice versa, in the form of a square with an odd side and a gnomon whose area is equal to the area of a square with an even side (Figure 2.2).

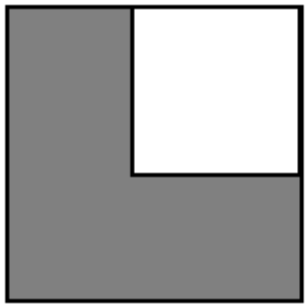

Figure 2.1

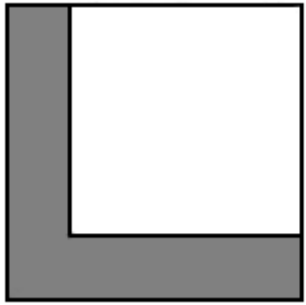

Figure 2.2

Two gnomons that make up one primitive Pythagorean triple together are called connected gnomons. Connected gnomons are constructed by mapping both squares to their own gnomons (Figure 2.3). At the same time, a larger gnomon absorbs a smaller gnomon.

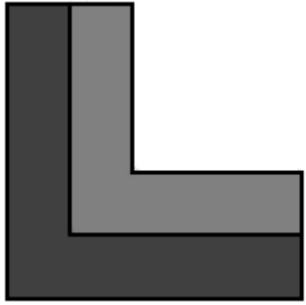

Figure 2.3

The parameters of the gnomon are its thickness **T** and the length of the side **L**.

A primitive Pythagorean triple is given

$$x^2 + y^2 = z^2$$

If $y < x$, then the thickness gnomon **T1** obtained from $x^2$ will be equal to

$$T1 = z - y = l^2$$

The thickness **T2** of the gnomon obtained from $y^2$ will be equal to

$$T2 = z - x = 2t^2$$

The lengths of the sides of both gnomons are equal to each other and equal in magnitude to **z**

$$L1 = L2 = z.$$

**Primitive Pythagorean triples and their representation via arithmetic progressions**

Let's imagine a primitive Pythagorean triple in the form of a square and a gnomon placed on it. Take a square with an odd side **x**. Then the area of the gnomon can be represented as the sum of an arithmetic progression with the first term **2x + 1**. Each subsequent member will be two units larger than the previous one. The number of such terms in the arithmetic progression is equal to the thickness of the gnomon

$$T2 = 2t^2.$$

Now take a square with an even side **y**. Then the area of the gnomon built on it can be represented as the sum of an arithmetic progression with the first term equal to $2y + 1$. Each subsequent member will be two units larger than the previous one. The number of such terms in the arithmetic progression is equal to the thickness of the gnomon

$$T1 = l^2.$$

For **y** < **x**, all the terms of the gnomon, and their number is **T2**, will be equal, respectively, to the last terms in the arithmetic progression representing $x^2$. And, conversely, for **x** < **y**, all the terms of the gnomon, and their number is **T1**, will be equal, respectively, to the last terms in the arithmetic progression representing $y^2$.

This representation in the form of an arithmetic progression of each gnomon fully corresponds to the picture of the gnomon absorbing a larger area of the associated gnomon of a smaller area.

**General Pythagorean triples**

Multiply all the elements of a primitive Pythagorean triple by an integer coefficient **k**.

Let's construct a primitive Pythagorean triple $(x, y, z)$. Let's draw it as a square with side **z**. Inside this square, in the upper right row, we place a square with the side **y**. Let's add the gnomon **U** to the inner square. The area of a gnomon is equal to the area of a square with side **x**. The thickness of the gnomon is $l^2$. We place the square $z^2$ in the cells of the square lattice with the side **kz** (Fig. 3).

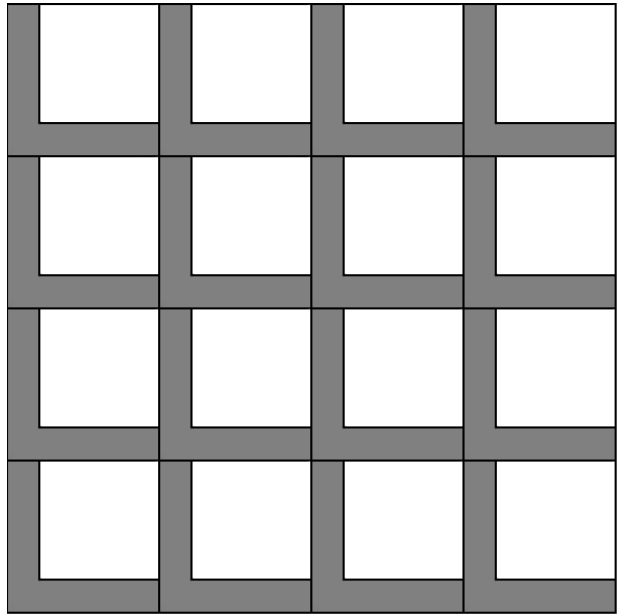

Figure 3 - Square lattice with side $kz$, when $k = 4$

Let's put together all the squares with the y side on the right, and on the left and at the bottom we will draw the total gnomon from the gnomons of each square lattice cell (Fig. 4).

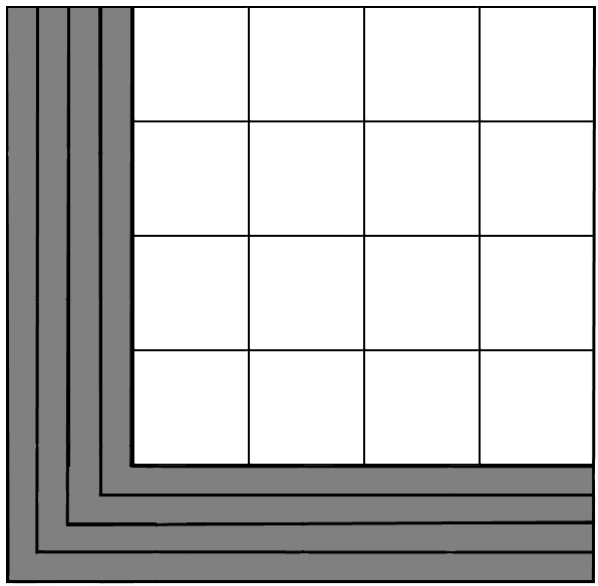

Figure 4 - Square lattice with side $kz$, when $k = 4$

We see that the thickness of the gnomon has become equal to $kl^2$. The area of the summing square $z^2$ has increased $k^2$ times. The area of the square with the side **y** has increased $k^2$ times. Consequently, the area of the gnomon also increased $k^2$ times. Similarly, the construction of a gnomon with an area of $y^2$. Its thickness will be $2kt^2$. The lengths of both gnomons will be equal to $kz$.

Thus, when multiplying all the elements of the Pythagorean triple by an integer factor $k$, the thickness of each gnomon in the corresponding constructions increases by $k$ times.

**Conclusion**

In this paper, an order relation is constructed on the set of primitive Pythagorean triples using two parameters. This allows for the solving of various tasks related to the ordered processing of these objects.

Formulas for the one-to-one mapping of the elements of primitive Pythagorean triples into the parameters of the partition of the side of the generating square are obtained. This allows for the immediate placing of the primitive Pythagorean triple in its place in the table.

A one-to-one correspondence is found between the elements of primitive Pythagorean triples on one side and the square plus gnomon pairs, as well as the pair of two connected gnomons.

A one-to-one correspondence between gnomon elements and arithmetic progressions is found. It is shown that the connectedness of gnomons is consistent with the overlap of the terms of arithmetic progressions representing gnomons.

In this paper, we propose a method for converting a primitive Pythagorean triple, represented as a square plus a gnomon, into a general Pythagorean triple, when multiplying all the elements of a primitive Pythagorean triple by a factor $k$. At the same time, for a general Pythagorean triple, the construction of a gnomon from individual gnomons of a primitive Pythagorean triple is geometrically shown.

**Fragment of a table of primitive Pythagorean triples constructed with increasing parameters $S, t(S)$**

| N.$n_i$ | S | t | l | x | y | z |
|---|---|---|---|---|---|---|
| 1.1 | 2 | 1 | 1 | 3 | 4 | 5 |
| 2.1 | 4 | 2 | 1 | 5 | 12 | 13 |
| 3.1 | 6 | 1 | 3 | 15 | 8 | 17 |
| 3.2 |   | 3 | 1 | 7 | 24 | 25 |
| 4.1 | 8 | 4 | 1 | 9 | 40 | 41 |
| 5.1 | 10 | 1 | 5 | 35 | 12 | 37 |
| 5.2 |   | 5 | 1 | 11 | 60 | 61 |
| 6.1 | 12 | 2 | 3 | 21 | 20 | 29 |
| 6.2 |   | 6 | 1 | 13 | 84 | 85 |
| 7.1 | 14 | 1 | 7 | 63 | 16 | 65 |
| 7.2 |   | 7 | 1 | 15 | 112 | 113 |
| 8.1 | 16 | 8 | 1 | 17 | 144 | 145 |
| 9.1 | 18 | 1 | 9 | 99 | 20 | 101 |
| 9.2 |   | 9 | 1 | 19 | 180 | 181 |
| 10.1 | 20 | 2 | 5 | 45 | 28 | 53 |
| 10.2 |   | 10 | 1 | 21 | 220 | 221 |
| 11.1 | 22 | 1 | 11 | 143 | 24 | 145 |
| 11.2 |   | 11 | 1 | 23 | 264 | 265 |
| 12.1 | 24 | 4 | 3 | 33 | 56 | 65 |
| 12.2 |   | 12 | 1 | 25 | 312 | 313 |
| 13.1 | 26 | 1 | 13 | 195 | 28 | 197 |
| 13.2 |   | 13 | 1 | 27 | 364 | 365 |
| 14.1 | 28 | 2 | 7 | 77 | 36 | 85 |
| 14.2 |   | 14 | 1 | 29 | 420 | 421 |
| 15.1 | 30 | 1 | 15 | 255 | 32 | 257 |
| 15.2 |   | 3 | 5 | 55 | 48 | 73 |
| 15.3 |   | 5 | 3 | 39 | 80 | 89 |
| 15.4 |   | 15 | 1 | 31 | 480 | 481 |
| 16.1 | 32 | 16 | 1 | 33 | 544 | 545 |
| 17.1 | 34 | 1 | 17 | 223 | 36 | 225 |
| 17.2 |   | 17 | 1 | 35 | 612 | 613 |
| 18.1 | 36 | 2 | 9 | 117 | 44 | 125 |
| 18.2 |   | 18 | 1 | 37 | 684 | 685 |
| 19.1 | 38 | 1 | 19 | 399 | 40 | 401 |
| 19.2 |   | 19 | 1 | 39 | 760 | 761 |
| 20.1 | 40 | 4 | 5 | 65 | 72 | 97 |
| 20.2 |   | 20 | 1 | 41 | 840 | 841 |
| 21.1 | 42 | 1 | 21 | 483 | 44 | 485 |
| 21.2 |   | 3 | 7 | 91 | 60 | 109 |

| N.$n_i$ | s | t | l | x | y | z |
|---|---|---|---|---|---|---|
| 21.3 | | 7 | 3 | 51 | 140 | 149 |
| 21.4 | | 21 | 1 | 43 | 924 | 925 |
| 22.1 | 44 | 2 | 11 | 165 | 52 | 173 |
| 22.2 | | 22 | 1 | 45 | 1012 | 1013 |
| 23.1 | 46 | 1 | 23 | 575 | 48 | 577 |
| 23.2 | | 23 | 1 | 47 | 1104 | 1105 |
| 24.1 | 48 | 8 | 3 | 57 | 176 | 185 |
| 24.2 | | 24 | 1 | 49 | 1200 | 1201 |
| 25.1 | 50 | 1 | 25 | 675 | 52 | 677 |
| 25.2 | | 25 | 1 | 51 | 1300 | 1301 |
| 26.1 | 52 | 2 | 13 | 221 | 60 | 229 |
| 26.2 | | 26 | 1 | 53 | 1404 | 1405 |
| 27.1 | 54 | 1 | 27 | 783 | 56 | 785 |
| 27.2 | | 27 | 1 | 55 | 1512 | 1513 |
| 28.1 | 56 | 4 | 7 | 105 | 88 | 137 |
| 28.2 | | 28 | 1 | 57 | 1624 | 1625 |
| 29.1 | 58 | 1 | 29 | 899 | 60 | 901 |
| 29.2 | | 29 | 1 | 59 | 1740 | 1741 |
| 30.1 | 60 | 2 | 15 | 285 | 68 | 293 |
| 30.2 | | 6 | 5 | 85 | 132 | 157 |
| 30.3 | | 10 | 3 | 69 | 260 | 269 |
| 30.4 | | 30 | 1 | 61 | 1860 | 1861 |
| 31.1 | 62 | 1 | 31 | 1023 | 64 | 1025 |
| 31.2 | | 31 | 1 | 63 | 1984 | 1985 |
| 32.1 | 64 | 32 | 1 | 65 | 2112 | 2113 |
| 33.1 | 66 | 1 | 33 | 1155 | 68 | 1157 |
| 33.2 | | 3 | 11 | 187 | 84 | 205 |
| 33.3 | | 11 | 3 | 75 | 308 | 317 |
| 33.4 | | 33 | 1 | 67 | 2244 | 2245 |
| 34.1 | 68 | 2 | 17 | 357 | 76 | 365 |
| 34.2 | | 34 | 1 | 69 | 2380 | 2381 |
| 35.1 | 70 | 1 | 35 | 1295 | 72 | 1297 |
| 35.2 | | 5 | 7 | 119 | 120 | 169 |
| 35.3 | | 7 | 5 | 95 | 168 | 193 |
| 35.4 | | 35 | 1 | 71 | 2520 | 2521 |
| 36.1 | 72 | 4 | 9 | 153 | 104 | 185 |
| 36.2 | | 36 | 1 | 73 | 2664 | 2665 |
| 37.1 | 74 | 1 | 37 | 1443 | 76 | 1445 |
| 37.2 | | 37 | 1 | 75 | 2812 | 2813 |
| 38.1 | 76 | 2 | 19 | 437 | 84 | 445 |
| 38.2 | | 38 | 1 | 77 | 2964 | 2965 |
| 39.1 | 78 | 1 | 39 | 1599 | 80 | 1601 |
| 39.2 | | 3 | 13 | 247 | 96 | 265 |
| 39.3 | | 13 | 3 | 87 | 416 | 425 |
| 39.4 | | 39 | 1 | 79 | 3120 | 3121 |
| 40.1 | 80 | 8 | 5 | 125 | 208 | 233 |

| N.$n_i$ | s | t | l | x | y | z |
|---|---|---|---|---|---|---|
| 40.2 |  | 40 | 1 | 81 | 3280 | 3281 |
| 41.1 | 82 | 1 | 41 | 1763 | 84 | 1765 |
| 41.2 |  | 41 | 1 | 83 | 3444 | 3445 |
| 42.1 | 84 | 2 | 21 | 525 | 92 | 533 |
| 42.2 |  | 6 | 7 | 133 | 156 | 205 |
| 42.3 |  | 14 | 3 | 93 | 476 | 485 |
| 42.4 |  | 42 | 1 | 85 | 3612 | 3613 |
| 43.1 | 86 | 1 | 43 | 1935 | 88 | 1937 |
| 43.2 |  | 43 | 1 | 87 | 3784 | 3785 |
| 44.1 | 88 | 4 | 11 | 209 | 120 | 241 |
| 44.2 |  | 44 | 1 | 89 | 3960 | 3961 |
| 45.1 | 90 | 1 | 45 | 2115 | 92 | 2117 |
| 45.2 |  | 5 | 9 | 171 | 140 | 221 |
| 45.3 |  | 9 | 5 | 115 | 252 | 277 |
| 45.4 |  | 45 | 1 | 91 | 4140 | 4141 |
| 46.1 | 92 | 2 | 23 | 621 | 100 | 629 |
| 46.2 |  | 46 | 1 | 93 | 4324 | 4325 |
| 47.1 | 94 | 1 | 47 | 2303 | 96 | 2305 |
| 47.2 |  | 47 | 1 | 95 | 4512 | 4513 |
| 48.1 | 96 | 16 | 3 | 105 | 608 | 617 |
| 48.2 |  | 48 | 1 | 97 | 4704 | 4705 |
| 49.1 | 98 | 1 | 49 | 2499 | 100 | 2501 |
| 49.2 |  | 49 | 1 | 99 | 4900 | 4901 |
| 50.1 | 100 | 2 | 25 | 725 | 108 | 733 |
| 50.2 |  | 50 | 1 | 101 | 5100 | 5101 |

# Example of mapping elements of a primitive Pythagorean triple through the parameters of the partition of the side of the generating square

In the archaeological collection of Columbia University, there is a cuneiform tablet dating from about 1500 BC. When it was examined, it turned out that it contains a list of Pythagorean triples. In this list, in addition to the triple (3,4,5), there is, for example, the triple (4961,6480,8161) [2]. This shows that the list was compiled by some method other than trial and error. Thus, the Babylonians had knowledge of the Pythagorean theorem a thousand years before Pythagoras and possessed methods of finding Pythagorean triples.

Let's map the elements of this triple (4961,6480,8161) to the parameters of the partition of the side of the generating square.

$$8161 - 6480 = 1681 = 41^2 = l^2;$$

$$l = 41;$$

$$4961 - 1681 = 3280 = 2 \cdot 41 \cdot t;$$

$$t = \frac{3280}{2 \cdot 41} = 40;$$

$$S = 2 \cdot t \cdot l = 2 \cdot 40 \cdot 41 = 3280;$$

So we have the side of the generating square $S = 3280$ with a partition $2tl = 2 \cdot 40 \cdot 41$.

Hence the primitive Pythagorean triple is easily constructed:

$$y = S + 2t^2 = 6480;$$

$$x = S + l^2 = 4961;$$

$$z = S + 2t^2 + l^2 = 8161.$$